\documentclass[reqno]{amsart}
\usepackage{amsmath, colortbl}
\usepackage{amsfonts}
\usepackage{amssymb}
\usepackage{mathrsfs}
\usepackage{graphicx}
\usepackage{flafter}
\newtheorem{thm}{Theorem}[section]
\newtheorem{cor}[thm]{Corollary}
\newtheorem{lem}[thm]{Lemma}

\theoremstyle{definition}

\newtheorem{rem}[thm]{Remark}
\numberwithin{equation}{section}

\begin{document}

\thispagestyle{empty}

\vspace{1 true cm} \large{
\title[ Regularity criterion for 3D N-S
Equations ] { Regularity criterion for  3D Navier-Stokes \\
Equations in  Besov spaces }%
\author[Daoyuan Fang \ Chenyin Qian]{Daoyuan Fang \hspace{0.5cm}Chenyin Qian$^*$\\
Department of Mathematics, Zhejiang University,\\ Hangzhou, 310027,
China}

\thanks{$^*$E-mail addresses: dyf@zju.edu.cn (D. Fang), qcyjcsx@163.com(C. Qian)}

\begin{abstract}

  Several regularity criterions  of Leray-Hopf weak solutions $u$ to the 3D Navier-Stokes
equations  are obtained. The results show that a weak solution $u$
becomes regular  if the gradient of  velocity component
$\nabla_{h}{u}$ (or $ \nabla{u_3}$) satisfies the additional
conditions in the  class of $L^{q}(0,T;
\dot{B}_{p,r}^{s}(\mathbb{R}^{3}))$, where
$\nabla_{h}=(\partial_{x_{1}},\partial_{x_{2}})$ is the horizontal
gradient operator. Besides, we also consider the anisotropic
regularity criterion for the weak solution of Navier-Stokes
equations in $\mathbb{R}^3$. Finally, we also get a further
regularity criterion, when give the sufficient condition on
$\partial_3u_3$.
\end{abstract}

\maketitle

\section{Introduction}
 In  the present  paper, we  address sufficient  conditions for the regularity of  weak solutions of the Cauchy problem for the
Navier-Stokes equations in $ \mathbb{R}^{3}\times (0,T)$:
\begin{equation} \label{a}
 \left\{\begin{array}{l}
\displaystyle \frac{\partial u}{\partial t}-\nu \Delta
u+(u\cdot\nabla)u+\nabla p=0,
\ \mbox{\ in}\ \mathbb{R}^{3}\times (0,T),\\
\displaystyle\nabla\cdot u=0,\hspace{3.52cm} \mbox{\ in}\ \mathbb{R}^{3}\times (0,T),\\
\displaystyle u(x, 0)=u_{0},\hspace{0.2cm}\mbox{\ in}\ \mathbb{R}^{3},\\
\end{array}
\right. \end{equation} where $u=(u_{1},u_{2},u_{3}):
\mathbb{R}^{3}\times (0,T)\rightarrow \mathbb{R}^{3}$ is the
velocity field, $ p: \mathbb{R}^{3}\times (0,T)\rightarrow
\mathbb{R}^{3}$ is a scalar pressure, and $u_{0}$ is the initial
velocity field, $\nu>0$ is the viscosity.  We set
$\nabla_{h}=(\partial_{x_{1}},\partial_{x_{2}})$ as the horizontal
gradient operator and
$\Delta_{h}=\partial_{x_{1}}^{2}+\partial_{x_{2}}^{2}$ as the
horizontal Laplacian, and $\Delta$ and $\nabla$ are the usual
Laplacian and the gradient operators, respectively.  Here we use the
classical notations
$$
(u\cdot\nabla)v=\sum_{i=1}^{3}u_{i}\partial_{x_{i}}v_{k}, \ (
k=1,2,3),\ \ \ \nabla\cdot u=\sum_{i=1}^{3}\partial_{x_{i}}u_{i},
$$
 and for sake of simplicity,  we denote $\partial_{x_{i}}$
by $\partial_{{i}}$.

\par   It is well known that the weak solution of the Navier-Stokes equations \eqref{a} is  unique  and
regular in two dimensions. However, in three dimensions, the
regularity problem of weak solutions of Navier-Stokes equations is
an outstanding open problem in mathematical fluid mechanics.  Strong
solutions are known to exist for a short interval of
time whose length depends on the initial data. Moreover, this strong
solution is known to be unique and to depend continuously on the
initial data (see, for example, \cite{[21]}, \cite{[22]}). Let us
recall the definition of Leray-Hopf weak solution. We set
$$
\mathcal {V}=\{\phi: \ \mbox{ the 3D vector valued}\ C_{0}^{\infty}
\ \mbox{functions and}\ \nabla\cdot\phi=0\},
$$
which will form the space of test functions. Let $H$ and $V$ be the
closure spaces of $\mathcal {V}$  under $L^2$-topology, and
 under $H^1$-topology, respectively.
\par
For $u_0\in H$,  the existence of weak solutions of \eqref{a} was
established by Leray  \cite{[14]} and Hopf in  \cite{[9]}, that is,
$u$
satisfies the following properties:\\
(i) $u\in C_{w}([0,T); H)\cap L^{2}(0,T; V)$, and $\partial_{t}u\in
L^{1}(0,T; V^{\prime})$, where $V^{\prime}$ is the dual space of
$V$;\\
(ii) $u$ verifies (1.1) in the sense of distribution, i.e., for
every test function $\phi\in C^{\infty}([0,T);\mathcal {V})$, and
for almost every $t, t_{0}\in (0,T)$, we have
$$
\begin{array}{ll}
 \displaystyle&\displaystyle\int_{\mathbb{R}^{3}} u(x,t)\cdot\phi(x,t)dx-\int_{\mathbb{R}^{3}}
u(x,t_{0})\cdot\phi(x,t_{0})dx\vspace{2mm}\\
\displaystyle &\ \ \ \
=\displaystyle\int_{t_{0}}^{t}\int_{\mathbb{R}_{3}}[u(x,t
)\cdot(\phi_{t}(x,t)+\nu\Delta\phi(x,t))]dxds\vspace{2mm}\\
&\ \ \ \ \ \ \
+\displaystyle\int_{t_{0}}^{t}\int_{\mathbb{R}_{3}}[(u(x,t)\cdot\nabla)\phi(x,t)]\cdot
u(x,t))]dxds
\end{array}
$$
 (iii) The energy inequality,
i.e.,
$$
\|u(\cdot,t)\|_{L^{2}}^{2}+2\nu\int_{t_0}^{t}\|\nabla
u(\cdot,s)\|_{L^{2}}^{2}ds\leq\|u_{0}\|_{L^{2}}^{2},
$$
for every $t$ and almost every $t_{0}$.

It is well known, if $u_{0}\in V$, a weak solution becomes  strong
solution of (1.1) on $(0, T)$ if, in addition, it satisfies
$$
u\in C([0,T); V)\cap L^{2}(0,T; H^{2}) \ \mbox{and}\
\partial_{t}u\in L^{2}(0,T; H).
$$

\par  Researchers are interested in the classical problem of finding sufficient conditions for
weak solutions of (1.1) such that the weak solutions become regular,
and the first result is usually referred as Prodi-Serrin conditions
(see \cite{[18]} and \cite{[20]}), which states that if a weak
solution $u$ is in the class of
$$ u\in L^{t}(0,T; L^{s}(\mathbb{R}^{3})),\ \
\frac{2}{t}+\frac{3}{s}=1,\ s\in[3,\infty],
$$ then the weak solution becomes regular.
\par
The full regularity of weak solutions can also be proved under
alternative assumptions on the gradient of the velocity $\nabla u$.
Specifically (see \cite{[1]}), if
\begin{equation}\label{zb}
\nabla u\in L^{t}(0,T; L^{s}(\mathbb{R}^{3})),\ \
\frac{2}{t}+\frac{3}{s}=2,\ s\in[\frac{3}{2},\infty].
\end{equation}
H. Beir$\tilde{\mbox{a}}$o da Veiga \cite{[28]} extended Serrin's
regularity criterion to the vorticity  $w=\mbox{curl}u$ showing that
if
\begin{equation}\label{b}
\omega\in L^{q}(0,T; L^r) \ \mbox{with} \ \frac{2}{q}+\frac{3}{r}=2,
\ \frac{3}{2}<r<\infty,
\end{equation}
then $u$ is  regular. In the marginal case $r=\infty$, H.
Kozono and Y. Taniuchi \cite{[29]} proved the regularity of weak
solutions under the condition
\begin{equation}\label{c}
\omega\in L^{1}(0,T; BMO),
\end{equation}
where $BMO$ is the space of bounded mean oscillation defined by
$$
f\in L^{1}_{{loc}}(\mathbb{R}^3), \
\sup_{x,R}\frac{1}{|B_{R}|}\int_{B_{R}(x)}|f(y)-\bar{f}_{B_{R}(x)}|dy<\infty,
$$
where $\bar{f}_{B_{R}(x)}$ is the average of $f$ over
$B_{R}(x)=\{y\in\mathbb{R}^3; |x-y|<R\}$. 
\par
Recently, the study of the regularity of weak solution involving
Besov space becomes popular.  For example, by establishing the
logarithmic Sobolev inequality in Besov spaces, H. Kozono, T. Ogawa
and Y. Taniuchi \cite{[32]} refined the  above two conditions  to
\begin{equation}\label{d}
\omega\in L^{q}(0,T; \dot{B}_{r,\infty}^{0})\ \mbox{with}\
\frac{2}{q}+\frac{3}{r}=2, \ {3}\leq r\leq\infty.
\end{equation}
Here and thereafter, $\dot{B}_{p,q}^{s}$ stands  for the homogeneous
Besov space, see Section 2 for the definition. On the other hand,
Chen and Zhang in \cite{[33]} proved  regularity criterion  by
imposing only the two-component vorticity field. More precisely,  they proved the regularity
of weak solutions in the class of
\begin{equation}\label{e}
\tilde{\omega}\in L^{q}(0,T; \dot{B}_{r,\sigma}^{0})\ \mbox{with}\ \tilde{\omega}= (\omega_1,\omega_2,0),
\frac{2}{q}+\frac{3}{r}=2, \ \frac{3}{2}<r\leq\infty,
\sigma\leq\frac{2}{3}.
\end{equation}
In \cite{[34]}, H. Kozono and Y. Yatsu  showed that if the Leray-Hopf
weak solution $u$ of \eqref{a} satisfies
\begin{equation}\label{f}
\tilde{\omega}\in L^{1}(0,T; BMO),
\end{equation}
then $u$ is  regular. More generally, B. Yuan  and  B. Zhang in
\cite{[98]}   prove the weak solution $u$ became regular if $\omega$
satisfies
\begin{equation*}
{\omega}\in L^{\frac{2}{2-\alpha}}(0,T;
\dot{B}_{\infty,\infty}^{-\alpha}(\mathbb{R}^{3}))\bigcap
L^{\frac{2}{1-\alpha}}(0,T;
\dot{B}_{\infty,\infty}^{-1-\alpha}(\mathbb{R}^{3})) \ \mbox{for}\
0<\alpha<1.
\end{equation*}
As to the endpoint case, S. Gala in \cite{[97]} showed that if
\begin{equation*}
{\omega}\in L^{2}(0,T;
\dot{B}_{\infty,\infty}^{-1}(\mathbb{R}^{3})),
\end{equation*}
then the solution $u$ was regular. \par We point out that H. Kozono,
T. Ogawa and Y. Taniuchi
 in \cite{[32]} (Theorem 3.5) also got the full regularity of weak solutions under alternative
assumptions on the velocity $u$. More precisely, if  $u$ satisfies
\begin{eqnarray}\label{g}
u\in L^{2}(0,T; \dot{B}_{\infty,\infty}^{0}(\mathbb{R}^{3})),
\end{eqnarray}
then the solution $u$ is regular. More generally,   A.
Cheskidov and R. Shvydkoy in \cite{[68]} proved that the solution
becames smooth if
\begin{eqnarray}\label{212}
u\in L^{r}((0, T); B^{\frac{2}{r}-1}_{\infty,\infty}),
\end{eqnarray}
for some $r\in (2,\infty)$, where $B^{s}_{p,q}$ stands  for the
nonhomogeneous Besov space (for detail see \cite{[38]}).

 Motivated
by the mentioned above, in this article, we consider assumptions on
the gradient of velocity  $\nabla u$ or the gradient of velocity
component
 $\nabla_h u$ (or $\nabla{u_3}$).
 \par Our main results  can be stated in the
following:
\begin{thm}\label{t1.2}
Let $u$ be a Leray-Hopf weak solution to the 3D Navier-Stokes
equations \eqref{a} with the initial value $u_{0}\in V$.  Suppose
the vorticity $\omega=curl u$ satisfies the condition
\begin{eqnarray}\label{ii} \nabla u \in L^{2}(0,T;
\dot{B}_{\infty,\infty}^{-1}(\mathbb{R}^{3})). \end{eqnarray}
 Then $u$ is
regular.
\end{thm}
\begin{rem}
  By the definition of Besov space and  Bernstein
inequality (see \eqref{a2} in section 2), we have
\begin{eqnarray}\label{z}C^{-1}\|
u\|_{\dot{B}_{\infty,\infty}^{0}}\leq\|\nabla
u\|_{\dot{B}_{\infty,\infty}^{-1}}\leq C\|
u\|_{\dot{B}_{\infty,\infty}^{0}}.\end{eqnarray} From \eqref{z}, it
is obvious that
 condition the \eqref{ii}  is  equivalent to   \eqref{g}. Therefore, Theorem \ref{t1.2} is easy  to get  from the result of H.
Kozono, T. Ogawa and Y. Taniuchi  (Theorem 3.5) in \cite{[32]}.
\end{rem}
\begin{thm}\label{t1.3} Let  $u_{0}$ and  $u$ be as   in
Theorem \ref{t1.2}.  Suppose that one of the following conditions is true: \\
 (i) the Leray-Hopf weak solution $u$ satisfies
\begin{eqnarray}\label{j} \nabla_h u\in L^{\frac{8}{3}}(0,T;
\dot{B}_{\infty,\infty}^{-1}(\mathbb{R}^{3})).
\end{eqnarray}
(ii)    $u$ satisfies  the following condition
\begin{eqnarray}\label{k} \nabla u_3\in L^{\frac{8}{5-2s}}(0,T;
\dot{B}_{\infty,\infty}^{-s}(\mathbb{R}^{3}))\ \mbox{with}\ 0<s<1.
\end{eqnarray}
 Then $u$ is
regular.
\end{thm}
\begin{cor}\label{c1.1} Suppose that $u_{0}\in V$, and
$u$ is a Leray-Hopf weak solution to the 3D Navier-Stokes equations
\eqref{a}. If $u$ satisfies the condition
\begin{eqnarray}\label{l} u\in L^{\frac{8}{3}}(0,T;
\dot{B}_{\infty,\infty}^{0}(\mathbb{R}^{3})),
\end{eqnarray}
or, for any (small) positive real number $\varepsilon$,  satisfies
\begin{eqnarray}\label{m} u_3\in L^{\frac{8}{5-2s}}(0,T;
\dot{B}_{\infty,\infty}^{-s+1}(\mathbb{R}^{3}))\ \mbox{with}\ 0<s<1.
\end{eqnarray}
 Then $u$ is
regular.
\end{cor}
\begin{rem}
Theorem \ref{t1.3} pays attention to the case of  the gradient of
velocity component, the first result of it proves
the regularity criterion by imposing only the two-component of the
gradient of velocity field, namely the horizontal gradient
components. By \eqref{z}, the Corollary\ref{c1.1} is easy to get
from Theorem \ref{t1.3}, and we see that the first result of
Corollary\ref{c1.1} is also a consequence  of \eqref{g} or Theorem
\ref{t1.2}. However, the regularity result in case of $\nabla_h u$
satisfies \eqref{j} is not a trivial corollary of Theorem \ref{t1.2}.
\end{rem}

In framework of the Lebesgue spaces, the regularity criterion problem has  been in-deep study  with the 
conditions in terms of one component $\nabla u_3$ ( for example, see
\cite{[26]}, \cite{[727]}) or one directional derivative
$\partial_3u$ ( for example, see \cite{[37]}, \cite {[827]}).
Because of the embedding \begin{equation}\label{9090}
L^{p}(\mathbb{R}^{3}))\hookrightarrow
\dot{B}_{r,q}^{\frac{3}{r}-\frac{3}{p}}(\mathbb{R}^{3}))\hookrightarrow\dot{B}_{\infty,\infty}^{-\frac{3}{p}}(\mathbb{R}^{3})),
\end{equation}
with $2\leq p<r,q\leq\infty$, a natural idea is to extend the
regularity criterion results to the framework of Besov spaces. To our knowledge, there are some results  in term of the
whole velocity or the vorticity. Motivated by the literature
\cite{[37]}, \cite {[827]}, we  want to consider  additional
condition only on $\partial_1u_3,\partial_2u_3,\partial_3u_3$
instead of $\nabla u$ in the more general spaces.  Our result
reads as:
\begin{thm}\label{t1.4} Let  $u_{0}$ and  $u$ be as   in
Theorem \ref{t1.2}.  Suppose that the additional  conditions of  $u$
 are satisfied
\begin{eqnarray}\label{p}\partial_3u_i\in L^{2}(0,T;
\dot{B}_{\infty,\infty}^{-1}(\mathbb{R}^3)), \  i= 1,2,
\end{eqnarray}
and
\begin{eqnarray}\label{q}\
\partial_3u_3\in L^{\frac{4}{-5s+2}}(0,T;
\dot{B}_{\infty,\infty}^{-s}(\mathbb{R}^3))\ \mbox{with}\ 0<
s<\frac{2}{5}.
\end{eqnarray}
 Then $u$ is
regular.
\end{thm}

\begin{thm}\label{t1.5} Let  $u_{0}$ and  $u$ be as   in
Theorem \ref{t1.2}.  Suppose that  $u$ the additional  condition
\begin{eqnarray}\label{qq}\
\partial_3u_3\in L^{\frac{24}{-29s+8}}(0,T;
\dot{B}_{\infty,\infty}^{-s}(\mathbb{R}^3))\ \mbox{with}\ 0<
s<\frac{8}{29}.
\end{eqnarray}
 Then $u$ is
regular.
\end{thm}

\begin{rem} By the embedding
\eqref{9090}, we know that the condition \eqref{p} is corresponding
to the endpoint case of the Prodi-Serrin conditions in the class of
$L^{q}(0,T; L^{p}(\mathbb{R}^{3}))$ with $p=3$ and $q=2$, which is
consistent with \eqref{zb}. While the condition \eqref{q} is in
consistent with \eqref{zb}, however, we  give a same range of $q$
with $2<q<\infty$ when $0<s<\frac{2}{5}$. Furthermore, if we provide
the sufficient condition, only in terms of $\partial_3u_3$, we shall
have a more strict condition, which is shown in Theorem \ref{t1.5}.
One can see the corresponding $q$ varies from $3$ to $\infty$.
\end{rem}

\par For the convenience, we recall the following version of
the three-dimensional  Sobolev and Ladyzhenskaya inequalities in the
whole space $\mathbb{R}^{3}$ (see, for example, \cite{[8]},
\cite{[11]}). There exists a positive constant $C$ such that
\begin{equation}\label{r}
\begin{array}{ll}
 \displaystyle
\|u\|_{r}&\displaystyle
\leq C \|u\|_{2}^{\frac{6-r}{2r}}\|\partial_{1}u\|_{2}^{\frac{r-2}{2r}}\|\partial_{2}u\|_{2}^{\frac{r-2}{2r}}\|\partial_{3}u\|_{2}^{\frac{r-2}{2r}}\\
&\leq C \|u\|_{2}^{\frac{6-r}{2r}}\|\nabla
u\|_{2}^{\frac{3(r-2)}{2r}},
\end{array}
\end{equation}
for every $u\in H^{1}(\mathbb{R}^{3})$ and every $r\in[2,6]$, where
$C$ is a constant depending only on $r$. Taking $\nabla$div on both
sides of \eqref{a} for smooth $(u; p)$, one can obtain
$$
-\Delta(\nabla
p)=\sum_{i,j}^{3}\partial_{i}\partial_{j}(\nabla(u_{i}u_{j})),
$$
therefore, the Calderon-Zygmund inequality in $\mathbb{R}^{3}$ (see
\cite{[30]})
\begin{equation}\label{s}
\|\nabla p\|_{q}\leq C\||\nabla u||u|\|_{q},
 1<q<\infty,
\end{equation}
holds, where $C$ is a positive constant depending only on $q$. And
there is another  estimates for the pressure
\begin{equation}\label{y} \| p\|_{q}\leq C\|u\|_{2q}^2, \
 1<q<\infty,
\end{equation}
Recall also that if div$u=0$ then the vorticity
$\omega=$curl$u$$=\nabla\times u$ has the following estimates (see
\cite{[37]}):
\begin{equation}\label{t}
C\|\omega\|_{q}\leq \|\nabla u\|_q\leq C(q)\|\omega\|_{q},\
1<q<\infty.
\end{equation}
Moreover,   if   div$u=0$, the expression
\begin{equation}\label{u}
\Delta u=\nabla(\mbox{div}u)-\nabla\times(\nabla\times u)
\end{equation}
can be reduced to
\begin{equation}\label{v}
\Delta u=-\nabla\times(\nabla\times u)=-\nabla\times\omega.
\end{equation}
On the other hand, note that div$\omega=0$, applying \eqref{t}, we
have
\begin{equation}\label{w}
C\|\nabla\times\omega\|_{q}\leq \|\nabla\omega\|_q\leq
C(q)\|\nabla\times\omega\|_{q},\ 1<q<\infty.
\end{equation}
Therefore, by \eqref{t} and \eqref{v}, we have
\begin{equation}\label{x}
C\|\Delta u\|_{q}\leq \|\nabla\omega\|_q\leq C(q)\|\Delta u\|_{q},\
1<q<\infty.
\end{equation}
\section{{ Preliminaries}} We begin this section with some
notations and Lemmas, which is useful for us to prove the main
results. In order to define Besov spaces, we first introduce the
Littlewood-Paley decomposition theory. Let $\mathcal
{S}(\mathbb{R}^3)$ be the Schwartz class of rapidly decreasing
function, given $f\in \mathcal {S}(\mathbb{R}^3)$, its Fourier
transformation $ \mathcal {F}f=\hat{f}$  is defined by
$$
\hat{f}(\xi)=\int_{\mathbb{R}^3}e^{-ix\cdot\xi}f(x)dx,
$$
and its inverse Fourier transform $ \mathcal {F}^{-1}f=\check{f}$ is
defined by
$$
\check{f}(x)=(2\pi)^{-3}\int_{\mathbb{R}^3}e^{ix\cdot\xi}f(\xi)d\xi.
$$
More generally, the Fourier transform of any $f\in\mathcal
{S}^\prime(\mathbb{R}^3)$, the space of tempered distributions, is
given by
$$
\langle \hat{f}, g\rangle=\langle {f}, \hat{g}\rangle,
$$
 for any $g\in\mathcal {S}(\mathbb{R}^3)$. The
Fourier transform is a bounded linear bijection from $\mathcal
{S}^\prime$ to $\mathcal {S}^\prime$ whose inverse is also bounded.
We fix the notation
$$
\mathcal {S}_h=\{\phi\in\mathcal {S},
\int_{\mathbb{R}^3}\phi(x)x^\gamma dx=0,
|\gamma|=0,1,2,\cdot\cdot\cdot\}.
$$
Its dual is given by
$$
\mathcal {S}^\prime_h=\mathcal {S}^\prime/\mathcal
{S}^{\bot}_h=\mathcal {S}^\prime/\mathcal {P},
$$
where $\mathcal {P}$ is the space of polynomial. In other words, two
distributions in $\mathcal {S}^\prime_h$ are identified as the same
if their difference is a polynomial.
 Let us choose two nonnegative radial
functions $\chi, \varphi \in\mathcal {S}(\mathbb{R}^3)$ supported in
$\mathfrak{B}=\{ \xi\in\mathbb{R}^3: |\xi|\leq 4/3\}$ and
$\mathfrak{C}=\{\xi\in\mathbb{R}^3: 3/4\leq|\xi|\leq 8/3\}$
respectively, such that
$$
\sum_{j\in\mathbb{Z}}\varphi(2^{-j}\xi)=1, \forall
\xi\in\mathbb{R}^3\backslash\{0\},
$$
and
$$
\chi(\xi)+\sum_{j\geq 0}\varphi(2^{-j}\xi)=1, \forall
\xi\in\mathbb{R}^3.
$$
Let $h=\mathcal {F}^{-1}\varphi$ and $\tilde{h}=\mathcal
{F}^{-1}\chi$, and then we define the homogeneous dyadic blocks
$\dot{\varDelta}_j$ and the homogeneous low-frequency cut-off
operator $\dot{S}_j$ as follows:
$$
\dot{\varDelta}_j
u=\varphi(2^{-j}D)u=2^{3j}\int_{\mathbb{R}^3}h(2^{j}y)u(x-y)dy,
$$
and
$$
\dot{S}_j
u=\chi(2^{-j}D)=2^{3j}\int_{\mathbb{R}^3}\tilde{h}(2^{j}y)u(x-y)dy.
$$
Informally, $\dot{\varDelta}_j$ is a frequency projection to the
annulus $\{|\xi|\sim 2^j\}$, while $\dot{S}_j$ is a frequency
projection to the ball $\{|\xi|\lesssim 2^j\}$. And one can easily
verify that $\dot{\varDelta}_j\dot{\varDelta}_k f=0$ if $|j-k|\geq
2$. Especially for any $f\in L^{2}(\mathbb{R}^3)$, we have  the
Littlewood-Paley decomposition:
$$
f=\sum_{j=-\infty}^{+\infty}\dot{\varDelta}_j f.
$$
We now give the definitions of Besov spaces. Let $s\in R$, $p, q \in
[1,\infty]$, the homogeneous Besov space
$\dot{B}_{p,q}^{s}(\mathbb{R}^3)$ is defined by  the full-dyadic
decomposition. We say that $f\in \dot{B}_{p,q}^{s}(\mathbb{R}^3)$,
if $f\in {S}^\prime_h$ and
$$
\sum_{j=-\infty}^{+\infty}(2^{js}\|\dot{\varDelta}_j
f\|_{L^{p}})^{q}<\infty,
$$
with the norm $$ \|f\|_{\dot{B}_{p,q}^{s}}=\left\{\begin{array}{l}
\displaystyle(\sum_{j=-\infty}^{+\infty}2^{qjs}\|\dot{\varDelta}_j
f\|_{L^{p}}^{q})^{\frac{1}{q}}, \ 1\leq q<\infty,\\
\displaystyle \sup_{j\in \mathbb{Z}} 2^{js}\|\dot{\varDelta}_j
f\|_{L^{p}}, \ q=\infty.
\end{array}
\right.$$ It is of interest to note that the homogeneous Besov space
$\dot{B}_{2,2}^{s}(\mathbb{R}^3)$ is equivalent to the homogeneous
Sobolev space $\dot{H}^{s}(\mathbb{R}^3)$. The following Bernstein
inequalities will be used in the next section.
\begin{lem}\label{l2.1} (see \cite{[38]})
Let $\mathfrak{B}$ be a ball and $\mathfrak{C}$ an annulus. A
constant $C$ exists such that for any nonnegative integer $k$, and
couple (p,q) in $[1,\infty]^{2}$ with  $1\leq p\leq q$, and any
function $u $ of $L^p(\mathbb{R}^d),$ we have
\begin{equation}\label{a1}
\verb"Supp" \ \hat{u}\subset \lambda\mathfrak{B}\Rightarrow
\sup_{|\alpha|=k}\|\partial^{\alpha}u\|_{L^{q}}\leq
C^{k+1}\lambda^{k+d(\frac{1}{p}-\frac{1}{q})}\|u\|_{L^{p}},
\end{equation}
\begin{equation}\label{a2}
\verb"Supp"\ \hat{u}\subset \lambda\mathfrak{C}\Rightarrow
C^{-k-1}\lambda^k\|u\|_{L^{p}}\leq\sup_{|\alpha|=k}\|\partial^{\alpha}u\|_{L^{p}}\leq
C^{k+1}\lambda^{k}\|u\|_{L^{p}}.
\end{equation}
\end{lem}
\begin{lem}\label{l2.2}(see \cite{[38]})
Let $1\leq q<p<\infty$ and $\alpha$ be a positive real number. A
constant $C$ exists such that
\begin{equation}\label{a6}
\|f\|_{L^{p}}\leq
C\|f\|_{\dot{B}_{\infty,\infty}^{-\alpha}}^{1-\theta}\|f\|^{\theta}_{\dot{B}^{\beta}_{q,q}},
 \ with \ \beta=\alpha\left(\frac{p}{q}-1\right)\ and\  \theta=\frac{q}{p}.
\end{equation}
\end{lem}
In particular, for $\beta=1, q=2$ and $p=4$, we get $\alpha=1$ and
\begin{equation}\label{a3}
\|f\|_{L^{4}}\leq
C\|f\|_{\dot{B}_{\infty,\infty}^{-1}}^{1/2}\|f\|^{1/2}_{\dot{H}^{1}},
\end{equation}
and further, if we give  the suitable values to  parameters $\beta,
q, p, \alpha$,  we  get other inequalities, for example
\begin{equation}\label{a4}
\|f\|_{L^{6}}\leq
C\|f\|_{\dot{B}_{\infty,\infty}^{-\frac{1}{2}}}^{2/3}\|f\|^{1/3}_{\dot{H}^{1}},\
\|f\|_{L^{3}}\leq
C\|f\|_{\dot{B}_{\infty,\infty}^{-2}}^{1/3}\|f\|^{2/3}_{\dot{H}^{1}}.
\end{equation}
\begin{lem} (see \cite{[38]}) A constant $C$ exists which satisfies the following
properties. If $s_1$ and $s_2$ are real numbers such that $s_1 <s_2$
and $\theta\in (0, 1))$,  for any $(p, r) \in [1,\infty]^2$  and any
$f \in\dot{B}_{p,r}^{s_1}\bigcap\dot{B}_{p,r}^{s_2}$, then we have

\begin{equation}\label{a5}
\|f\|_{\dot{B}_{p,r}^{\theta s_1+(1-\theta)s_2}}\leq
C\|f\|_{\dot{B}_{p,r}^{ s_1}}^{\theta}\|f\|_{\dot{B}_{p,r}^{
s_2}}^{(1-\theta)}.
\end{equation}
\end{lem}
\section{{Proof of Main Results}}

In this section, under the assumptions of the Theorem \ref{t1.2},
Theorem \ref{t1.3} or Theorem \ref{t1.4} in Section 1 respectively,
we prove our main results.  First of all, we note that, by the
energy inequality, for Leray-Hopf weak solutions, we have (see, for
example, \cite{[21]}, \cite{[22]} for detail)
\begin{equation}\label{1}
\|u(\cdot,t)\|_{L^{2}}^{2}+2\nu\int_{0}^{t}\|\nabla
u(\cdot,s)\|_{L^{2}}^{2}ds\leq K_{1},\end{equation} for all $0<
t<T,$ where $K_{1}=\|u_{0}\|_{L^{2}}^{2}.$
\par  It is  well
known that there exists a unique strong solution  $u$ local in
time  if $u_{0}\in V$. In addition, this strong solution
$u\in C([0,T^{*});V)\cap L^{2}(0,T^{*}; H^{2}(\mathbb{R}^{3}))$ is
the only weak solution with the initial datum $u_{0}$, where
$(0,T^{*})$ is the maximal interval of existence of the unique
strong solution. If $T^{*}\geq T,$ then there is nothing to prove.
If, on the other hand, $T^{*}< T,$ then our strategy is to show that
the $H^{1}$ norm of this strong solution is bounded uniformly in
time over the interval $(0,T^{*})$, provided  additional conditions
in Theorem \ref{t1.2},  Theorem \ref{t1.3} or Theorem \ref{t1.4} in
Section 1 are valid. As a result the interval $(0,T^{*})$ can not be
a maximal interval of existence, and consequently $T^{*}\geq T,$
which concludes our proof. \par In order to prove the $H^{1}$ norm
of the strong solution $u$ is bounded on interval $(0,T^{*})$,
combing with the energy equality \eqref{1}, it is sufficient to
prove
\begin{equation} \label{2}
\|\nabla
u\|_{L^2}^{2}+\displaystyle{\nu}\displaystyle\int_{0}^{t}\|\Delta
u\|_{L^2}^{2}d\tau\leq C,\ \forall \ t\in(0, T^{*})
\end{equation}
where the constant $C$ depends on $T$, $K_{1} $.
\par

\textbf{Proof of Theorem \ref{t1.2}}  Taking the curl on \eqref{a},
we obtain
$$
\frac{\partial \omega}{\partial t}-\nu \Delta
\omega+(u\cdot\nabla)\omega-(\omega\cdot\nabla)u=0.
$$
We taking the inner product of above inequality with $\omega$  in
$L^{2}(\mathbb{R}^3)$, and by using of the H$\ddot{\mbox{o}}$lder's
and Young's inequalities, as well as \eqref{t}, \eqref{x} and
\eqref{a3}, we obtain
\begin{equation}\label{43}
\begin{array}{ll}
\displaystyle \frac{1}{2}\frac{d}{dt}\|
\omega\|_{L^2}^{2}+\nu\|\nabla \omega\|_{L^2}^{2}&=\displaystyle
\int_{\mathbb{R}^{3}}(\omega\cdot\nabla)u\cdot\omega dx\displaystyle
\\&=\displaystyle
\sum_{i,j=1}^{3}\int_{\mathbb{R}^{3}}\omega_i\partial_iu_j\omega_j dx\displaystyle \\
\displaystyle &\leq C\displaystyle\|\omega\|_{L^4}^2\|\nabla u\|_{L^2}\displaystyle\vspace{2mm}\\
\displaystyle &\leq C\displaystyle\|\nabla u\|_{L^4}^2\|\nabla u\|_{L^2}\displaystyle\vspace{2mm}\\
\displaystyle &\leq  C\|\nabla u\|_{\dot{B}_{\infty,\infty}^{-1}}\|\Delta u\|_{{L^2}}\|\nabla u\|_{L^2} \vspace{2mm}\\
\displaystyle &\leq  C\|\nabla u\|_{\dot{B}_{\infty,\infty}^{-1}}^2\|\nabla u\|_{L^2}^2+\displaystyle\frac{\nu}{2}\|\Delta u\|_{{L^2}}^2. \\
&\leq C\|\nabla u\|_{\dot{B}_{\infty,\infty}^{-1}}^2\|\omega\|_{L^2}^2+\displaystyle\frac{\nu}{2}\|\nabla \omega\|_{{L^2}}^2. \\
\end{array}
\end{equation}
Absorbing the first  term in right hand side and integrating the
above inequality, we obtain
\begin{equation}\label{44}
\|\omega\|_{L^2}^{2}+\nu\int_{0}^{t}\|\nabla
\omega\|_{L^2}^{2}\leq\|\omega_{0}\|_{L^2}^2+C\int_{0}^{t} \|\nabla
u\|_{\dot{B}_{\infty,\infty}^{-1}}^2\|\omega\|_{L^2}^2d\tau.
\end{equation}
Therefore, by Gronwall's inequality, one has
$$
\|\omega\|_{{L^2}}^{2}+\displaystyle{\nu}\displaystyle\int_{0}^{t}\|\nabla
\omega \|_{{L^2}}^{2}d\tau \leq
\|\omega_{0}\|_{L^2}^2\exp\left(C\int_{0}^{t}\|\nabla
u\|_{\dot{B}_{\infty,\infty}^{-1}}^2d\tau\right).
$$
 By using of Gronwall's inequality and
condition \eqref{ii}, we have $$ \omega\in L^{\infty}(0,T^{*};
L^{2}(\mathbb{R}^3))\cap L^{2}(0,T^{*}; H^{1}(\mathbb{R}^3)).$$
Therefore, by \eqref{t} and \eqref{x}, we get the $ H^1$ norm of the
strong solution $u$ is bounded on the maximal interval of existence
$(0, T^{*})$. This completes the proof of Theorem \ref{t1.2}.\\
 \textbf{Proof of Theorem
\ref{t1.3}} Firstly, we deal with (i). Taking the inner product of
the equation \eqref{a} with $-\Delta_{h}u$ in $L^{2}(\mathbb{R}^3)$,
By integrating by parts a few times and using the incompressibility
condition, we obtain
\begin{equation}\label{6}
\begin{array}{ll}
 \displaystyle \frac{1}{2}\frac{d}{dt}\|\nabla_{h}u\|_{L^2}^{2}+\nu\|\nabla_{h}\nabla u\|_{L^2}^{2}&=\displaystyle
 \int_{\mathbb{R}^{3}}[(u\cdot \nabla)u]\cdot\Delta_{h}u dx\\
 &\displaystyle
 =-\int_{\mathbb{R}^3}\sum_{k,j=1}^3\sum_{i=1}^{2}\partial_ju_k\partial_iu_j\partial_iu_kdx\displaystyle\\
 &\displaystyle\leq C \int_{\mathbb{R}^3}|\nabla u||\nabla_h
 u|^{2}dx.
\end{array}\end{equation}
Applying H$\ddot{\mbox{o}}$lder's  and Young's inequalities, as well
as \eqref{a3}, we have
\begin{equation}\label{7}
\begin{array}{ll}
 \displaystyle \frac{1}{2}\frac{d}{dt}\|\nabla_{h}u\|_{L^2}^{2}+\nu\|\nabla_{h}\nabla u\|_{L^2}^{2}&
 \leq\displaystyle
 \|\nabla_h
 u\|_{L^4}^{2}\|\nabla u\|_{L^2}\\
 &\displaystyle
 \leq \|\nabla_h
 u\|_{\dot{B}_{\infty,\infty}^{-1}}\|\nabla_h\nabla u\|_{L^2}\|\nabla u\|_{L^2}\displaystyle\\
 &\displaystyle
 \leq \|\nabla_h
 u\|_{\dot{B}_{\infty,\infty}^{-1}}^{2}\|\nabla u\|_{L^2}^2+\frac{\nu}{2}\|\nabla_h\nabla u\|_{L^2}^2.\displaystyle\\
\end{array}\end{equation}
Absorbing the first  term in right hand side and integrating the above
inequality, we obtain
\begin{equation}\label{8}
\|\nabla_{h}u\|_{L^2}^{2}+\displaystyle\nu\int_{0}^{t}\|\nabla_{h}\nabla
u\|_{L^2}^{2}d\tau\leq \displaystyle\|\nabla_{h}
u(0)\|_{L^2}^{2}+\int_{0}^{t} \|\nabla_h
 u\|_{\dot{B}_{\infty,\infty}^{-1}}^{2}\|\nabla u\|_{L^2}^2d\tau.
\end{equation}
Next, we  also use $-\Delta u$ as test function, and get
$$\begin{array}{ll} &\displaystyle \frac{1}{2}\frac{d}{dt}\|\nabla
u\|_{L^2}^{2}+\nu\|\Delta u\|_{L^2}^{2}\\&=\displaystyle
\sum_{i,j,k=1}^{3}\int_{\mathbb{R}^{3}}u_{i}\partial_{i}u_{j}\partial_{kk}u_{j} dx\displaystyle\\
&\displaystyle=\sum_{j=1}^{3}\int_{\mathbb{R}^3}u_3\partial_3u_j\Delta_{h}u_jdx+
\sum_{i=1}^{2}\sum_{j=1}^{3}\int_{\mathbb{R}^3}u_i\partial_iu_j\Delta
u_jdx+\sum_{j=1}^{3}\int_{\mathbb{R}^3}u_3\partial_3u_j\partial_{33}u_jdx\\
&=L_{1}(t)+L_{2}(t)+L_{3}(t)
\end{array}$$
The calculation has been shown  in  \cite{[26]}, for the convenience
of  readers, we list it
 below. By integrating by parts a few times and using the incompressibility
condition, we get $L_{1}(t), L_{2}(t), L_{3}(t)$ as follows
$$
\begin{array}{ll}\displaystyle L_{1}(t)&\displaystyle=-\sum_{j=1}^{3}\sum_{k=1}^{2}\int_{\mathbb{R}^3}\partial_ku_3\partial_3u_j\partial_ku_jdx-
\sum_{j=1}^{3}\sum_{k=1}^{2}\int_{\mathbb{R}^3}u_3\partial_{3k}u_j\partial_ku_jdx\\
&\displaystyle=-\sum_{j=1}^{3}\sum_{k=1}^{2}\int_{\mathbb{R}^3}\partial_ku_3\partial_3u_j\partial_ku_jdx+\frac{1}{2}
\sum_{j=1}^{3}\sum_{k=1}^{2}\int_{\mathbb{R}^3}\partial_3u_3\partial_{k}u_j\partial_ku_jdx,
\end{array}$$
$$
\begin{array}{ll}\
\displaystyle
L_{2}(t)&\displaystyle=-\sum_{i=1}^{2}\sum_{j=1}^{3}\sum_{k=1}^{3}\int_{\mathbb{R}^3}\partial_ku_i\partial_iu_j\partial_ku_jdx-
\sum_{i=1}^{2}\sum_{j=1}^{3}\sum_{k=1}^{3}\int_{\mathbb{R}^3}u_i\partial_{ik}u_j\partial_ku_jdx\\
&=\displaystyle-\sum_{i=1}^{2}\sum_{j=1}^{3}\sum_{k=1}^{3}\int_{\mathbb{R}^3}\partial_ku_i\partial_iu_j\partial_ku_jdx+\frac{1}{2}
\sum_{i=1}^{2}\sum_{j=1}^{3}\sum_{k=1}^{3}\int_{\mathbb{R}^3}\partial_iu_i\partial_{k}u_j\partial_ku_jdx,
\end{array}$$
$$
\begin{array}{ll}\
\displaystyle
L_{3}(t)&\displaystyle=-\frac{1}{2}\sum_{j=1}^{3}\int_{\mathbb{R}^3}\partial_3u_3\partial_3u_j\partial_3u_jdx=
\frac{1}{2}\sum_{j=1}^{3}\int_{\mathbb{R}^3}(\partial_1u_1+\partial_2u_2)\partial_3u_j\partial_3u_jdx.
\end{array}$$
Therefore,  by \eqref{r} and H$\ddot{\mbox{o}}$lder's inequalities,
for every $i\ (i=1,2,3)$ we have
\begin{equation} \label{9}\begin{array}{ll}\displaystyle
|L_{i}(t)|& \displaystyle\leq C\int_{\mathbb{R}^3}|\nabla_h u| |\nabla u|^{2}dx\vspace{1mm}\\
& \displaystyle\leq C\|\nabla_{h}u\|_{L^2}\|\nabla u\|_{L^4}^{2}\vspace{1mm}\\
&\displaystyle\leq C\|\nabla_{h}u\|_{L^2}\|\nabla
u\|_{L^2}^{\frac{1}{2}}\|\nabla_{h}\nabla u\|_{L^2}\|\Delta
u\|_{L^2}^{\frac{1}{2}},
\end{array}\end{equation}
and hance we have
\begin{equation} \label{10}\begin{array}{ll}\displaystyle \frac{1}{2}\frac{d}{dt}\|\nabla
u\|_{L^2}^{2}+\nu\|\Delta u\|_{L^2}^{2}\leq \displaystyle
C\|\nabla_{h}u\|_{L^2}\|\nabla
u\|_{L^2}^{\frac{1}{2}}\|\nabla_{h}\nabla u\|_{L^2}\|\Delta
u\|_{L^2}^{\frac{1}{2}}.
\end{array}\end{equation}
Integrating \eqref{10}, applying  H$\ddot{\mbox{o}}$lder's
inequality and combing \eqref{8} and \eqref{9}, we obtain
\begin{equation} \label{11}\begin{array}{ll}
&\|\nabla
u\|_{L^2}^{2}+\displaystyle2\nu\displaystyle\int_{0}^{t}\|\Delta
u\|_{L^2}^{2}d\tau\\
\displaystyle &\hspace{0.3cm} \leq\|\nabla
u(0)\|_{L^2}^{2}+\left(\sup_{0\leq s\leq t}\|\nabla_{h}
u\|_{L^2}\right)\left(\displaystyle\int_{0}^{t}\|\nabla
u\|_{L^2}^{2}d\tau \right)^{\frac{1}{4}}\\
&\ \ \ \ \ \times\left(\displaystyle\int_{0}^{t}\|\nabla_{h}\nabla
u\|_{L^2}^{2}d\tau\right)^{\frac{1}{2}}\left(\displaystyle\int_{0}^{t}\|\Delta
u\|_{L^2}^{2}d\tau\right)^{\frac{1}{4}}\\
\displaystyle &\hspace{0.3cm} \leq\|\nabla
u(0)\|_{L^2}^{2}+C\displaystyle\left(\int_{0}^{t}\|\nabla_h
 u\|_{\dot{B}_{\infty,\infty}^{-1}}^{2}\|\nabla u\|_{L^2}^2d\tau\right)\times
 \left(\int_{0}^{t}\|\Delta
 u\|_{L^2}^{2}d\tau\right)^{{\frac{1}{4}}}\\
 &\hspace{0.5cm} + \displaystyle\|\nabla_{h}
u(0)\|_{L^2}^{2}\left(\displaystyle\int_{0}^{t}\|\Delta
u\|_{L^2}^{2}d\tau\right)^{\frac{1}{4}}.
\end{array}\end{equation}
By using of the H$\ddot{\mbox{o}}$lder's and Young's inequalities,
it follows that
\begin{equation} \label{12}\begin{array}{ll}
&\|\nabla
u\|_{L^2}^{2}+\displaystyle{\nu}\displaystyle\int_{0}^{t}\|\Delta
u\|_{L^2}^{2}d\tau\\
\displaystyle &\hspace{1.3cm} \leq\|\nabla u(0)\|_{L^2}^{2}+
\displaystyle C\|\nabla_{h} u(0)\|_{L^2}^{8/3}\vspace{1mm}\\
& \hspace{1.5cm}+C\displaystyle\left(\int_{0}^{t}\|\nabla_h
 u\|_{\dot{B}_{\infty,\infty}^{-1}}^{2}\|\nabla u\|_{L^2}^2d\tau\right)^{4/3}\vspace{1mm}\\
\displaystyle &\hspace{1.3cm} \leq\|\nabla u(0)\|_{L^2}^{2}+
\displaystyle C\|\nabla_{h} u(0)\|_{L^2}^{8/3}\vspace{1mm}\\
& \hspace{1.5cm}+C\displaystyle\left(\int_{0}^{t}\|\nabla_h
 u\|_{\dot{B}_{\infty,\infty}^{-1}}^{\frac{8}{3}}\|\nabla u\|_{L^2}^2d\tau\right)\left
 (\int_{0}^{t}\|\nabla u\|_{L^{2}}^2d\tau\right)^{1/4}.
\end{array}\end{equation}
Thanks again to the energy inequality, we get
\begin{equation} \label{13}\begin{array}{ll}
\|\nabla
u\|_{L^2}^{2}+\displaystyle{\nu}\displaystyle\int_{0}^{t}\|\Delta
u\|_{L^2}^{2}d\tau\displaystyle&\leq\|\nabla u(0)\|_{L^2}^{2}+
\displaystyle C\|\nabla_{h} u(0)\|_{L^2}^{8/3}\vspace{1mm}\\
& \ \ \ \ \ +C\displaystyle\int_{0}^{t}\|\nabla_h
 u\|_{\dot{B}_{\infty,\infty}^{-1}}^{\frac{8}{3}}\|\nabla u\|_{L^2}^2d\tau.
\end{array}\end{equation}
Therefore, by using of Gronwall's inequality, we finally obtain
$$\begin{array}{ll}
&\|\nabla
u\|_{L^2}^{2}+\displaystyle{\nu}\displaystyle\int_{0}^{t}\|\Delta
u\|_{L^2}^{2}d\tau\\&  \  \ \ \ \displaystyle\leq \left(\|\nabla
u(0)\|_{L^2}^{2}+ \displaystyle C\|\nabla_{h}
u(0)\|_{L^2}^{8/3}\right)\exp\left(C\int_{0}^{t}\|\nabla_h
 u\|_{\dot{B}_{\infty,\infty}^{-1}}^{\frac{8}{3}}d\tau\right),
\end{array}$$
by condition \eqref{j}, we get the $ H^1$ norm of the strong
solution $u$ is bounded on the maximal interval of existence $(0,
T^{*})$. This completes the proof of (i).  \par Now we prove (ii).
Taking the inner product of the equation \eqref{a} with
$-\Delta_{h}u$ in $L^{2}(\mathbb{R}^3)$, we have (see  \cite{[2]}
for detail)
\begin{equation}\label{14}
\begin{array}{ll}
 \displaystyle \frac{1}{2}\frac{d}{dt}\|\nabla_{h}u\|_{L^2}^{2}+\nu\|\nabla_{h}\nabla u\|_{L^2}^{2}&=\displaystyle
 \int_{\mathbb{R}^{3}}[(u\cdot \nabla)u]\cdot\Delta_{h}u dx\\
&\leq\displaystyle
 C\int_{\mathbb{R}^{3}}|u_{3}||\nabla u||\nabla_{h}\nabla u|dx
\end{array}\end{equation}
By using  of the Littlewood-Paley decomposition, we decompose $u_3$
as follows:
\begin{equation}\label{16}
u_3=\sum_{j=-\infty}^{+\infty}\dot{\varDelta}_{j}u_3=\sum_{j<[\sigma]}\dot{\varDelta}_{j}u_3+\sum_{j\geq
[\sigma]+1}\dot{\varDelta}_{j}u_3,
\end{equation}
where $\sigma$ is a  real number determined later, and $[\cdot]$
denotes the integer part of $\sigma$. Therefore, \eqref{14} becomes
\begin{equation}\label{17}
\displaystyle
\frac{1}{2}\frac{d}{dt}\|\nabla_{h}u\|_{L^2}^{2}+\nu\|\nabla_{h}\nabla
u\|_{L^2}^{2}\leq I_{1}(t)+I_{2}(t),
\end{equation}
with
$$
I_{1}(t)=C\sum_{j<[\sigma]}\int_{\mathbb{R}^{3}}|\dot{\varDelta}_{j}u_3||\nabla
u||\nabla_{h}\nabla u|dx,
$$
$$
I_{2}(t)=C\sum_{j\geq
[\sigma]+1}\int_{\mathbb{R}^{3}}|\dot{\varDelta}_{j}u_3||\nabla
u||\nabla_{h}\nabla u|dx.
$$
In what following, we estimate $I_{1}(t)$ and $I_{2}(t)$. For
$I_{1}(t)$, by using of the H$\ddot{\mbox{o}}$lder's and Young's
inequalities, as well as Lemma 2.1, we have
\begin{equation}\label{18}
\begin{array}{ll}
 \displaystyle I_{1}(t)& \displaystyle\leq C\sum_{j<[\sigma]}
\|\dot{\varDelta}_{j}u_3\|_{L^{\infty}}\|\nabla
u\|_{L^2}\|\nabla_{h}\nabla u\|_{L^{2}}\\
& \displaystyle \leq C
\left(\sum_{j<[\sigma]}2^{\frac{3}{2}j}\right)
\|u_3\|_{L^{2}}\|\nabla u\|_{L^2}\|\nabla_{h}\nabla u\|_{L^{2}}\vspace{2mm}\\
&\displaystyle \leq C 2^{\frac{3}{2}[\sigma]}\|u_3\|_{L^{2}}\|\nabla
u\|_{L^2}\|\nabla_{h}\nabla u\|_{L^{2}}.\vspace{2mm}\\
&\displaystyle \leq C 2^{3\sigma}\|u_3\|_{L^{2}}^2\|\nabla
u\|_{L^2}^2+ \frac{\nu}{4}\|\nabla_{h}\nabla u\|_{L^{2}}^2,
\end{array}\end{equation}
the last inequality, we use the fact that $[\sigma]\le\sigma$.  As
to $I_{2}(t)$, we take the same strategy to $I_{1}(t)$, by the
definition of norm of the Besov space, for any $0<\varepsilon<1,$ we
have
\begin{equation}\label{19}
\begin{array}{ll}
 \displaystyle I_{2}(t)& \displaystyle\leq C\sum_{j\geq [\sigma]+1}
\|\dot{\varDelta}_{j}u_3\|_{L^{\infty}}\|\nabla
u\|_{L^2}\|\nabla_{h}\nabla u\|_{L^{2}}\\
& \displaystyle \leq C \sum_{j\geq [\sigma]+1}2^{-\varepsilon
j}2^{(-1+\varepsilon)j}\|\dot{\varDelta}_{j}\nabla
u_3\|_{L^{\infty}}\|\nabla u\|_{L^2}\|\nabla_{h}\nabla u\|_{L^{2}}\\
& \displaystyle \leq C \left(\sum_{j\geq [\sigma]+1}2^{-\varepsilon
j}\right)\|\nabla
u_3\|_{\dot{B}_{\infty,\infty}^{-1+\varepsilon}}\|\nabla u\|_{L^2}\|\nabla_{h}\nabla u\|_{L^{2}}\\
&\displaystyle \leq C2^{-2\varepsilon ([\sigma]+1)} \|\nabla
u_3\|_{\dot{B}_{\infty,\infty}^{-1+\varepsilon}}^2\|\nabla
u\|_{L^2}^2+ \frac{\nu}{4}\|\nabla_{h}\nabla u\|_{L^{2}}^2.
\\
&\displaystyle \leq C2^{-2\varepsilon \sigma} \|\nabla
u_3\|_{\dot{B}_{\infty,\infty}^{-1+\varepsilon}}^2\|\nabla
u\|_{L^2}^2+ \frac{\nu}{4}\|\nabla_{h}\nabla u\|_{L^{2}}^2,
\end{array}\end{equation}
the last inequality, we use the fact that $\sigma<[\sigma]+1$.
Inserting \eqref{18} and \eqref{19} into \eqref{14} to obtain
\begin{equation}\label{0000}\displaystyle
\frac{d}{dt}\|\nabla_{h}u\|_{L^2}^{2}+\nu\|\nabla_{h}\nabla
u\|_{L^2}^{2}\leq C 2^{-2\varepsilon \sigma}\|\nabla
u_3\|_{\dot{B}_{\infty,\infty}^{-1+\varepsilon}}^2\|\nabla
u\|_{L^2}^2+ C2^{3\sigma} \|u_3\|_{L^{2}}^2\|\nabla u\|_{L^2}^2,
\end{equation}
Now, we choose $\sigma$ such that
$$2^{\frac{3}{2}\sigma}\|u_3\|_{L^2}\|\nabla u\|_{L^2}
=2^{-\varepsilon \sigma}\|\nabla
u_3\|_{B^{-1+\varepsilon}_{\infty,\infty}}\|\nabla u\|_{L^2},$$ then
we have $$C2^{-\varepsilon \sigma}\|\nabla
u_3\|_{B^{-1+\varepsilon}_{\infty,\infty}}\|\nabla u\|_{L^2}\leq
C\|u_3\|_{L^2}^{\frac{\varepsilon}{\frac{3}{2}+\varepsilon}}
 \|\nabla u_3\|_{B^{-1+\varepsilon}_{\infty,\infty}}^{\frac{\frac{3}{2}}{\frac{3}{2}+\varepsilon}}\|\nabla u\|_{L^2}.$$
 Integrating \eqref{0000}, combing  above  two inequalities and the energy
inequality \eqref{1} we have
\begin{equation}\label{20}
\|\nabla_{h}u\|_{{L^2}}^{2}+\displaystyle\nu\int_{0}^{t}\|\nabla_{h}\nabla
u\|_{{L^2}}^{2}d\tau\leq \displaystyle\|\nabla_{h}
u(0)\|_{{L^2}}^{2}+C+C\int_{0}^{t}\|\nabla
u_3\|_{\dot{B}_{\infty,\infty}^{-1+\varepsilon}}^{\frac{3}{\frac{3}{2}+\varepsilon}}\|\nabla
u\|_{L^2}^2d\tau.
\end{equation}
Integrating \eqref{10}, applying  H$\ddot{\mbox{o}}$lder's
inequality and combing \eqref{1} and \eqref{20}, we obtain
\begin{equation} \label{21}\begin{array}{ll}
&\|\nabla
u\|_{L^2}^{2}+\displaystyle2\nu\displaystyle\int_{0}^{t}\|\Delta
u\|_{L^2}^{2}d\tau\\
\displaystyle &\hspace{0.3cm} \leq\|\nabla
u(0)\|_{L^2}^{2}+\left(\sup_{0\leq s\leq t}\|\nabla_{h}
u\|_{L^2}\right)\left(\displaystyle\int_{0}^{t}\|\nabla
u\|_{L^2}^{2}d\tau\right)^{\frac{1}{4}}\\&\ \ \ \ \
\times\left(\displaystyle\int_{0}^{t}\|\nabla_{h}\nabla
u\|_{L^2}^{2}d\tau\right)^{\frac{1}{2}}\left(\displaystyle\int_{0}^{t}\|\Delta
u\|_{L^2}^{2}d\tau\right)^{\frac{1}{4}}\\
\displaystyle &\hspace{0.3cm} \leq\|\nabla
u(0)\|_{L^2}^{2}+C\displaystyle\left(\int_{0}^{t}\|\nabla
u_3\|_{\dot{B}_{\infty,\infty}^{-1+\varepsilon}}^{\frac{3}{\frac{3}{2}+\varepsilon}}\|\nabla
u\|_{L^2}^2d\tau\right)\times \left(\int_{0}^{t}\|\Delta
 u\|_{L^2}^{2}d\tau\right)^{{\frac{1}{4}}}\\
 &\hspace{0.5cm} + \displaystyle\left(\|\nabla_{h}
u(0)\|_{L^2}^{2}+C\right)\left(\displaystyle\int_{0}^{t}\|\Delta
u\|_{L^2}^{2}d\tau\right)^{\frac{1}{4}}.
\end{array}\end{equation}
By using of the H$\ddot{\mbox{o}}$lder's and Young's inequalities,
it follows that
\begin{equation} \label{22}\begin{array}{ll}
&\|\nabla
u\|_{L^2}^{2}+\displaystyle{\nu}\displaystyle\int_{0}^{t}\|\Delta
u\|_{L^2}^{2}d\tau\\
\displaystyle &\hspace{1.3cm} \leq\|\nabla u(0)\|_{L^2}^{2}+
\displaystyle C\left(\|\nabla_{h} u(0)\|_{L^2}^{8/3}+1\right)\vspace{1mm}\\
& \hspace{1.5cm}+C\displaystyle\left(\int_{0}^{t}\|\nabla
u_3\|_{\dot{B}_{\infty,\infty}^{-1+\varepsilon}}^{\frac{3}{\frac{3}{2}+\varepsilon}}\|\nabla
u\|_{L^2}^2d\tau\right)^{4/3}\vspace{1mm}\\
\displaystyle &\hspace{1.3cm} \leq\|\nabla u(0)\|_{L^2}^{2}+
\displaystyle C\left(\|\nabla_{h} u(0)\|_{L^2}^{8/3}+1\right)\vspace{1mm}\\
& \hspace{1.5cm}+C\displaystyle\left(\int_{0}^{t}\|\nabla
u_3\|_{\dot{B}_{\infty,\infty}^{-1+\varepsilon}}^{{\frac{4}{\frac{3}{2}+\varepsilon}}}\|\nabla
u\|_{L^2}^2d\tau\right)\left(\int_{0}^{t}\|\nabla
u\|_{L^{2}}^2d\tau\right)^{1/4}.
\end{array}\end{equation}
Thanks again to the energy inequality, we get
\begin{equation} \label{23}\begin{array}{ll}
\|\nabla
u\|_{L^2}^{2}+\displaystyle{\nu}\displaystyle\int_{0}^{t}\|\Delta
u\|_{L^2}^{2}d\tau\displaystyle&\leq\|\nabla u(0)\|_{L^2}^{2}+
\displaystyle C\left(\|\nabla_{h} u(0)\|_{L^2}^{8/3}+1\right)\vspace{1mm}\\
& \ \ \ \ \ +C\displaystyle\int_{0}^{t}\|\nabla
u_3\|_{\dot{B}_{\infty,\infty}^{-1+\varepsilon}}^{{\frac{4}{\frac{3}{2}+\varepsilon}}}\|\nabla
u\|_{L^2}^2d\tau.
\end{array}\end{equation}
If we set $s=1-\varepsilon$, then we have
$$
{\frac{4}{\frac{3}{2}+\varepsilon}}=\frac{8}{5-2s}\ \mbox{with}\
0<s<1.
$$
Therefore, by using of Gronwall's inequality, we finally obtain
$$\begin{array}{ll}
&\|\nabla
u\|_{L^2}^{2}+\displaystyle{\nu}\displaystyle\int_{0}^{t}\|\Delta
u\|_{L^2}^{2}d\tau\\&  \  \ \ \ \displaystyle\leq \left(\|\nabla
u(0)\|_{L^2}^{2}+ \displaystyle C(\|\nabla_{h}
u(0)\|_{L^2}^{8/3}+1)\right)\exp\left(C\int_{0}^{t}\|\nabla
u_3\|_{\dot{B}_{\infty,\infty}^{-s}}^{\frac{8}{5-2s}}d\tau\right),
\end{array}$$
by condition \eqref{k}, we get the $ H^1$ norm of the strong
solution $u$ is bounded on the maximal interval of existence $(0,
T^{*})$. This completes the proof of Theorem \ref{t1.3}.
\par \textbf{Proof
of Theorem \ref{t1.4}} We also split the proof into two steps.
Recalling that in Theorem 1.2, the vorticity $\omega$ satisfies
$$
\begin{array}{ll}
\displaystyle \frac{1}{2}\frac{d}{dt}\|
\omega\|_{L^2}^{2}+\nu\|\nabla \omega\|_{L^2}^{2}&=\displaystyle
\int_{\mathbb{R}^{3}}(\omega\cdot\nabla)u\cdot\omega
dx\displaystyle=\displaystyle
\sum_{i,j=1}^{3}\int_{\mathbb{R}^{3}}\omega_i\partial_iu_j\omega_j dx.\displaystyle \\
\end{array}$$
Since $\omega=(\omega_1, \omega_2,
\omega_3)=(\partial_2u_3-\partial_3u_2,
\partial_3u_1-\partial_1u_3, \partial_1u_2-\partial_2u_1),$
we put the detail computation to get
\begin{equation} \label{24}\begin{array}{ll}
\displaystyle &\displaystyle\frac{1}{2}\frac{d}{dt}\|
\omega\|_{L^2}^{2}+\nu\|\nabla \omega\|_{L^2}^{2}\vspace{1mm}\\
&\ \ \ \ \ \ \ \ \displaystyle
=\int_{\mathbb{R}^{3}}\partial_3u_2\partial_3u_2\partial_1u_1dx
-\int_{\mathbb{R}^{3}}\partial_1u_2\partial_3u_2\partial_3u_1dx\displaystyle\vspace{1mm}\\
&\ \ \ \ \ \ \ \ \ \  \  \displaystyle
+\int_{\mathbb{R}^{3}}\partial_3u_1\partial_3u_1\partial_2u_2dx
-\int_{\mathbb{R}^{3}}\partial_2u_1\partial_3u_2\partial_3u_1dx\displaystyle\vspace{1mm}\\
&\ \ \ \ \ \ \ \ \ \ \ \displaystyle-\int_{\mathbb{R}^{3}}u_3\{
\partial_2(\partial_2u_3\partial_1u_1-2\partial_3u_2\partial_1u_1
+\partial_1u_2\partial_3u_1\vspace{1mm}\\
&\ \ \ \ \ \ \ \ \ \ \  \ \ \ \ \
\displaystyle-\partial_1u_2\partial_1u_3\displaystyle-\partial_2u_1\partial_1u_3
+\partial_1u_1\partial_3u_1)\}dx \vspace{1mm}\\
&\ \ \ \ \ \ \ \ \ \  \ \displaystyle-\int_{\mathbb{R}^{3}}u_3\{
\partial_1(\partial_3u_2\partial_1u_2+\partial_1u_3\partial_2u_1
-2\partial_2u_2\partial_3u_1\vspace{1mm}\\
&\ \ \ \ \ \ \ \ \ \ \ \ \ \ \ \ \
\displaystyle+\partial_2u_2\partial_1u_3-\partial_3u_2\partial_1u_3\displaystyle+\partial_2u_1\partial_3u_2
)\}dx \\
&\ \ \ \ \ \ \ \ \ \ \ \ \displaystyle
-\int_{\mathbb{R}^{3}}\sum_{i=1}^{3}
\omega_iu_3\partial_i\omega_3dx,
\end{array}\end{equation}
then, by Young's inequality, we have
\begin{equation} \label{25}\begin{array}{ll}
\displaystyle &\displaystyle\frac{1}{2}\frac{d}{dt}\|
\omega\|_{L^2}^{2}+\nu\|\nabla \omega\|_{L^2}^{2}\vspace{1mm}\\
&\ \ \ \ \ \ \ \ \leq \displaystyle
\int_{\mathbb{R}^{3}}|\partial_3u_2|^2|\partial_1u_1|dx+
\int_{\mathbb{R}^{3}}|\partial_3u_2||\partial_1u_2||\partial_3u_1|dx\vspace{1mm}\\
&\ \ \ \ \ \ \ \ \ \ \ \  + \displaystyle
\int_{\mathbb{R}^{3}}|\partial_3u_1|^2|\partial_2u_2|dx+
\int_{\mathbb{R}^{3}}|\partial_3u_2||\partial_2u_1||\partial_3u_1|dx\vspace{1mm}\\
&\ \ \ \ \ \ \ \ \ \ \ \  +\displaystyle C
\int_{\mathbb{R}^{3}}|u_3||\nabla
u||\Delta u|dx.\\
&\ \ \ \ \ \ \ \ \leq\displaystyle
C\int_{\mathbb{R}^{3}}|\partial_3u_2|^2|\nabla
u|dx+C\int_{\mathbb{R}^{3}}|\partial_3u_1|^2|\nabla u|dx\\
&\ \ \ \ \ \ \ \ \ \ \ \  +\displaystyle C
\int_{\mathbb{R}^{3}}|u_3||\nabla
u||\Delta u|dx.\\
&\ \ \ \ \ \ \ \ =K_{1}(t)+K_{2}(t)+K_{3}(t).\\
\end{array}\end{equation}
Next, we estimate $K_i(t)$ one by one, $i=1,2,3$. Applying
H$\ddot{\mbox{o}}$lder's and Young's inequalities, as well as
\eqref{t}, \eqref{x} and \eqref{a3}, we have
\begin{equation} \label{26}\begin{array}{ll}
K_{1}(t)& \displaystyle =C\int_{\mathbb{R}^{3}}|\partial_3u_2|^2|\nabla u|dx\vspace{1mm}\\
&\displaystyle\leq C\|\partial_3u_2\|_{L^4}^2\|\nabla u\|_{L^2}\vspace{1mm}\\
&\displaystyle\leq
C\|\partial_3u_2\|_{\dot{B}_{\infty,\infty}^{-1}}\|\nabla\partial_3u_2\|_{L^2}\|\nabla
u\|_{L^2}\vspace{1mm}\\
&\displaystyle\leq
C\|\partial_3u_2\|_{\dot{B}_{\infty,\infty}^{-1}}^{2}\|\nabla
u\|_{L^2}^2+\frac{\nu}{8}\|\Delta u\|_{{L^2}}^{2}\vspace{1mm}\\
&\displaystyle \leq
C\|\partial_3u_2\|_{\dot{B}_{\infty,\infty}^{-1}}^{2}\|\omega\|_{L^2}^2
+\frac{\nu}{8}\|\nabla \omega\|_{L^2}^{2},\vspace{1mm}
\end{array}\end{equation}
\begin{equation} \label{27}\begin{array}{ll}
K_{2}(t)\displaystyle
=C\int_{\mathbb{R}^{3}}|\partial_3u_1|^2|\nabla u|dx\displaystyle
\leq
C\|\partial_3u_1\|_{\dot{B}_{\infty,\infty}^{-1}}^{2}\|\omega\|_{L^2}^2
+\frac{\nu}{8}\|\nabla \omega\|_{{L^2}}^{2}.\vspace{1mm}
\end{array}\end{equation}
Now, we estimate $K_{3}(t)$, again, applying
H$\ddot{\mbox{o}}$lder's and Young's inequalities, as well as
\eqref{r} and \eqref{x}, we obtain
\begin{equation} \label{28}\begin{array}{ll}
K_{3}(t)& \displaystyle =\displaystyle C
\int_{\mathbb{R}^{3}}|u_3||\nabla
u||\Delta u|dx\vspace{1mm}\\
&\displaystyle\leq C\|u_3\|_{L^q}\|\nabla u\|_{L^p}\|\Delta u\|_{L^2}\vspace{1mm}\\
&\displaystyle\leq C \|u_3\|_{L^q}\|\nabla
u\|_{{L^2}}^{\frac{6-p}{2p}}\|\Delta
u\|_{{L^2}}^{\frac{5p-6}{2p}}\vspace{1mm}\\
&\displaystyle\leq C \|u_3\|_{L^q}^{\frac{4p}{6-p}}\|\nabla
u\|_{{L^2}}^{2}+\frac{\nu}{16}\|\Delta u\|_{{L^2}}^{2}\vspace{1mm}\\
&\displaystyle= C \|u_3\|_{L^q}^{\frac{2q}{q-3}}\|\nabla
u\|_{{L^2}}^{2}+\frac{\nu}{16}\|\nabla\omega\|_{{L^2}}^{2},\vspace{1mm}
\end{array}\end{equation}
where $p$ and $q$ satisfy
$$
\frac{1}{p}+\frac{1}{q}=\frac{1}{2} \ \mbox{with}\ 2\leq p<6, q>3.
$$
For $\|u_3\|_{L^q}$, we  have the following estimate. We use
$|u_{3}|^{q-2}u_{3}$ as test function  in the equation \eqref{a} for
$u_{3}.$ By using of Gagliardo-Nirenberg and
H$\ddot{\mbox{o}}$lder's inequalities, and applying  the inequality
\eqref{y}, we have
\begin{equation}\label{29}
\begin{array}{ll}
 \displaystyle \frac{1}{q}\frac{d}{dt}\|u_{3}\|_{L^q}^{q}+C\nu\|\nabla|u_{3}|^{{\frac{q}{2}}}\|_{L^2}^{2}
 &=\displaystyle
 -\int_{\mathbb{R}^{3}}\partial_{3}p|u_{3}|^{q-2}u_{3}dx\displaystyle  \vspace{1mm}\\
 &\leq C\displaystyle\int_{\mathbb{R}^{3}}|p||u_{3}|^{q-2}|\partial_{3}u_{3}|dx\displaystyle  \vspace{2mm}\\
 \displaystyle &\leq\displaystyle
 C \|p\|_{L^\mu}\|u_{3}\|_{L^{q}}^{q-2}\|\partial_{3}u_{3}\|_{{L^\beta}}\ \displaystyle \vspace{2mm}\\
 \displaystyle &\leq\displaystyle   C \|u
 \|_{L^{2\mu}}^2\|u_{3}\|_{L^q}^{q-2}\|\partial_3u_{3}\|_{L^\beta}\displaystyle\vspace{2mm}\\
 \displaystyle &\leq\displaystyle C\| u\|_{{L^2}}^{\frac{3-\mu}{\mu}}\|\nabla
 u\|_{L^2}^{\frac{3(\mu-1)}{\mu}}
 \|u_{3}\|_{L^q}^{q-2}\|\partial_{3}u_{3}\|_{L^\beta}.
\end{array}
\end{equation}
The parameters in above inequality satisfy
\begin{equation}\label{30}
\frac{1}{\mu}+\frac{1}{\beta}+\frac{q-2}{q}=1, \ \mbox{with} \
1\leq\mu\leq3,\  \beta>7.
\end{equation}
By \eqref{x}, \eqref{a6} and the energy inequality, The inequality
\eqref{29} implies that
\begin{equation}\label{31}
\begin{array}{ll}
\displaystyle
\frac{1}{2}\frac{d}{dt}\|u_{3}\|_{L^q}^{2}\displaystyle
&\leq\displaystyle C\| u\|_{L^2}^{\frac{3-\mu}{\mu}}\|\nabla
 u\|_{L^2}^{\frac{3(\mu-1)}{\mu}}\|\partial_{3}u_{3}\|_{L^\beta}\\
\displaystyle& \leq\displaystyle C\|
u\|_{L^2}^{\frac{3-\mu}{\mu}}\|\nabla
 u\|_{L^2}^{\frac{3(\mu-1)}{\mu}}\|\partial_{3}u_{3}\|_{\dot{B}_{\infty,\infty}^{\frac{2}{2-\beta}}}^{1-\frac{2}{\beta}}
 \|\nabla\partial_{3}u_{3}\|_{L^2}^{\frac{2}{\beta}}\\
 \displaystyle& \leq\displaystyle C\|
u\|_{L^2}^{\frac{3-\mu}{\mu}}\|\nabla
u\|_{L^2}^{\frac{3(\mu-1)}{\mu}}\|\partial_{3}u_{3}\|_{\dot{B}_{\infty,\infty}^{\frac{2}{2-\beta}}}^{1-\frac{2}{\beta}}
 \|\Delta u\|_{L^2}^{\frac{2}{\beta}},\end{array}
\end{equation} Integrating \eqref{31}, applying \eqref{1} and
H$\ddot{\mbox{o}}$lder's inequality, we have
\begin{equation}\label{32}
\begin{array}{ll}
\displaystyle \|u_{3}\|_{L^q}^{2}\displaystyle
&\leq\displaystyle\|u_{3}(0)\|_{L^q}^2+ C\int_{0}^{t}\|
u\|_{L^2}^{\frac{3-\mu}{\mu}}\|\nabla
u\|_{L^2}^{\frac{3(\mu-1)}{\mu}}\|\partial_{3}u_{3}\|_{\dot{B}_{\infty,\infty}^{\frac{2}{2-\beta}}}^{1-\frac{2}{\beta}}
 \|\Delta u\|_{L^2}^{\frac{2}{\beta}}d\tau\vspace{1mm}\\
 &\leq\displaystyle C\left(\int_{0}^{t}\|\partial_{3}u_{3}\|_{\dot{B}_{\infty,\infty}^{\frac{2}{2-\beta}}}^{(1-\frac{2}{\beta})
 (\frac{\beta}{\beta-1})}\|\nabla
 u\|_{L^2}^{\frac{3\beta(\mu-1)}{\mu(\beta-1)}}d\tau\right)^{\frac{\beta-1}{\beta}}\\
 &\displaystyle \ \ \ \  \ \ \times\left(\int_{0}^{t}\|\Delta u\|
 |_{L^2}^{2}d\tau\right)^{\frac{1}{\beta}}+\|u_{3}(0)\|_{L^q}^2\\
\end{array}
\end{equation}
We choose
\begin{equation} \label{33}\mu=\frac{3\beta}{\beta+2},
\end{equation}
then
$$
\frac{3\beta(\mu-1)}{\mu(\beta-1)}=2.
$$
Combing \eqref{28} and \eqref{32}, by energy inequality, and using
H$\ddot{\mbox{o}}$lder's  and Young's inequalities, we get
\begin{equation} \label{34}\begin{array}{ll}
\displaystyle\int_{0}^{t}K_{3}(\tau)d\tau &\displaystyle\leq  C
\left(\sup_{0\leq\tau\leq
t}\|u_3\|_{L^q}^2\right)^{\frac{q}{q-3}}\int_0^t\|\nabla
u\|_{{L^2}}^{2}d\tau+\frac{\nu}{16}\int_0^t\|\nabla\omega\|_{{L^2}}^{2}d\tau\\
&\displaystyle\leq
C\left(\int_{0}^{t}\|\partial_{3}u_{3}\|_{\dot{B}_{\infty,\infty}^{\frac{2}{2-\beta}}}^{(1-\frac{2}{\beta})
 (\frac{\beta}{\beta-1})}
\|\nabla u\|_{L^2}^{2}d\tau\right)^{\frac{q(\beta-1)}{\beta(q-3)}}
\left(\int_{0}^{t}\|\nabla\omega\|_{L^2}^{2}d\tau\right)^{\frac{q}{\beta(q-3)}}\\
& \ \ \ \ \ \displaystyle +
C\|u_{3}(0)\|_{L^q}^{\frac{2q}{q-3}}+\displaystyle\frac{\nu}{16}\displaystyle\int_0^t\|\nabla\omega\|_{{L^2}}^{2}d\tau\\
&\displaystyle\leq
C\left(\int_{0}^{t}\|\partial_{3}u_{3}\|_{\dot{B}_{\infty,\infty}^{\frac{2}{2-\beta}}}^{(1-\frac{2}{\beta})
 (\frac{\beta}{\beta-1})}
\|\nabla
u\|_{L^2}^{2}d\tau\right)^{\frac{q(\beta-1)}{\beta(q-3)-q}}\displaystyle+C\|u_{3}(0)\|_{L^q}^{\frac{2q}{q-3}}\\
&\ \ \ \ \ +\displaystyle\frac{\nu}{8}\displaystyle\int_0^t\|\nabla\omega\|_{{L^2}}^{2}d\tau\\
&\displaystyle\leq
C\left(\int_{0}^{t}\|\partial_{3}u_{3}\|_{\dot{B}_{\infty,\infty}^{\frac{2}{2-\beta}}}^{\frac{q(\beta-2)}{\beta(q-3)-q}}
\|\nabla u\|_{L^2}^{2}d\tau\right)\left(\int_{0}^{t}\|\nabla u\|_{L^2}^{2}d\tau\right)^{\frac{3\beta}{q(\beta-1)}}\\
& \ \ \ \ \ \displaystyle +
C\|u_{3}(0)\|_{L^q}^{\frac{2q}{q-3}}+\displaystyle\frac{\nu}{8}\displaystyle\int_0^t\|\nabla\omega\|_{{L^2}}^{2}d\tau\\
&\displaystyle\leq\displaystyle
C\int_{0}^{t}\|\partial_{3}u_{3}\|_{\dot{B}_{\infty,\infty}^{\frac{2}{2-\beta}}}^{\frac{q(\beta-1)}{\beta(q-3)-q}}
\|\omega\|_{L^2}^{2}d\tau+C\|u_{3}(0)\|_{L^q}^{\frac{2q}{q-3}}\\ & \
\ \ \ \ \displaystyle +
\frac{\nu}{8}\int_{0}^{t}\|\nabla\omega\|_{L^2}^{2}d\tau,\vspace{1mm}\\
\end{array}\end{equation}
in above inequality, we note that $q$ and $\beta$ satisfy the
condition
$$
\frac{q}{\beta(q-3)}<1.
$$

Integrating \eqref{25},  combing \eqref{26}, \eqref{27}, \eqref{34},
and applying Young's inequality, we have
\begin{equation} \label{35}\begin{array}{ll}
\displaystyle &\displaystyle\|
\omega\|_{L^2}^{2}+\nu\int_{0}^{t}\|\nabla \omega\|_{L^2}^{2}d\tau\vspace{1mm}\\
&\ \ \ \ \ \ \ \ \displaystyle\leq
C\int_{0}^{t}\|\partial_3u_2\|_{\dot{B}_{\infty,\infty}^{-1}}^{2}\|\omega\|_{L^2}^2d\tau
+C\int_{0}^{t}\|\partial_3u_1\|_{\dot{B}_{\infty,\infty}^{-1}}^{2}\|\omega\|_{L^2}^2d\tau\\
&\ \ \ \ \ \ \ \ \displaystyle\ \ \ \
+C\int_{0}^{t}\|\partial_{3}u_{3}\|_{\dot{B}_{\infty,\infty}^{\frac{2}{2-\beta}}}^{\frac{q(\beta-1)}{\beta(q-3)-q}}
\|\omega\|_{L^2}^{2}d\tau\displaystyle+C\|u_{3}(0)\|_{L^q}^{\frac{2q}{q-3}}+\|\omega(0)\|_{L^2}^2.\\
\end{array}\end{equation}
If we denote $s=\frac{2}{\beta-2}$, then by \eqref{30} and
\eqref{33}, we have
$$
\frac{q(\beta-1)}{\beta(q-3)-q}=\frac{4}{-5s+2}\ \mbox{with}\
0<s<\frac{2}{5}
$$
By  using of  Gronwall's inequality, we obtain
\begin{equation}\label{36}\begin{array}{ll} \displaystyle &\displaystyle\|
\omega\|_{L^2}^{2}+\nu\int_{0}^{t}\|\nabla \omega\|_{L^2}^{2}d\tau\vspace{1mm}\\
&\ \ \ \ \ \ \ \ \displaystyle\leq
\left(C\|u_{3}(0)\|_{L^q}^{\frac{2q}{q-3}}+\|\omega(0)\|_{L^2}^2\right)
\exp\left(C\int_{0}^{t}\|\partial_{3}u_{3}\|_{\dot{B}_{\infty,\infty}^{-s}}^{\frac{4}{-5s+2}}
d\tau\right),\\
&\ \ \ \ \ \ \ \ \ \ \ \
\displaystyle\times\exp\left(C\int_{0}^{t}\|\partial_3u_1\|_{\dot{B}_{\infty,\infty}^{-1}}^{2}
d\tau\right)
\exp\left(C\int_{0}^{t}\|\partial_3u_2\|_{\dot{B}_{\infty,\infty}^{-1}}^{2}d\tau\right)
\end{array}\end{equation}
by the condition  \eqref{p} and \eqref{q},  we have $$ \omega\in
L^{\infty}(0,T^{*}; L^{2}(\mathbb{R}^3))\cap L^{2}(0,T^{*};
H^{1}(\mathbb{R}^3)).$$ Therefore, the $ H^1$ norm of the strong
solution $u$ is bounded on the maximal interval of existence $(0,
T^{*})$. This completes the proof of Theorem \ref{t1.4}.

\par \textbf{Proof
of Theorem \ref{t1.5}}
 Firstly, we begin with \eqref{14}, and by \eqref{r}, we have
\begin{equation}\label{101}
\begin{array}{ll}
 \displaystyle \frac{1}{2}\frac{d}{dt}\|\nabla_{h}u\|_{L^2}^{2}&+\nu\|\nabla_{h}\nabla
 u\|_{L^2}^{2}\\
&\leq\displaystyle
 C\int_{\mathbb{R}^{3}}|u_{3}||\nabla u||\nabla_{h}\nabla u|dx\\
&\displaystyle\leq C\|u_3\|_{L^q}\|\nabla u\|_{L^p}\|\nabla_h\nabla u\|_{L^2}\vspace{1mm}\\
&\displaystyle\leq C \|u_3\|_{L^q}\|\nabla
u\|_{{L^2}}^{\frac{6-p}{2p}}\|\Delta
u\|_{{L^2}}^{\frac{p-2}{2p}}\|\nabla_h\nabla u\|_{L^2}^{\frac{2(p-1)}{p}}\vspace{1mm}\\
&\displaystyle\leq C \|u_3\|_{L^q}^{p}\|\nabla
u\|_{{L^2}}^{\frac{6-p}{2}}\|\Delta u\|_{{L^2}}^{\frac{p-2}{2}}+\frac{\nu}{2}\|\nabla_h\nabla u\|_{L^2}^2,\vspace{1mm}\\
\end{array}\end{equation}
where $p$ and $q$ satisfy
$$
\frac{1}{p}+\frac{1}{q}=\frac{1}{2}\ \mbox{with}\ 2\leq p<5,
q>\frac{10}{3}.
$$
Therefore, integrating above inequality and using
H$\ddot{\mbox{o}}$lder's inequality, it follows that
\begin{equation}\label{102}
\begin{array}{ll}
 \displaystyle \|\nabla_{h}u\|_{L^2}^{2}&+\nu\displaystyle\int_0^t\|\nabla_{h}\nabla
 u\|_{L^2}^{2}\\
&\leq\displaystyle
 \displaystyle\|\nabla_{h}
u(0)\|_{{L^2}}^{2}+\displaystyle\int_0^t\|u_3\|_{L^q}^{p}\|\nabla
u\|_{{L^2}}^{\frac{6-p}{2}}\|\Delta u\|_{{L^2}}^{\frac{p-2}{2}}d\tau\vspace{1mm}\\
&\leq\displaystyle
 \displaystyle\|\nabla_{h}
u(0)\|_{{L^2}}^{2}+\displaystyle\left(\int_0^t\|u_3\|_{L^q}^{\frac{4p}{6-p}}\|\nabla
u\|_{{L^2}}^{2}d\tau\right)^{\frac{6-p}{4}}\displaystyle\left(\int_0^t\|\Delta u\|_{{L^2}}^{\frac{p-2}{2}}d\tau\right)^{\frac{p-2}{4}}\vspace{1mm}\\
\end{array}\end{equation}
For $\|u_3\|_{L^q}$, we  have the same estimate to \eqref{32}, in
which the parameters  satisfy \eqref{33} and
\begin{equation}\label{105}
\frac{1}{\mu}+\frac{1}{\beta}+\frac{q-2}{q}=1, \ \mbox{with} \
1\leq\mu\leq3,\  \beta>\frac{37}{4}.
\end{equation}

Next, in view of \eqref{10}, integrating \eqref{10}, applying
H$\ddot{\mbox{o}}$lder's  and Young's inequalities and combing
\eqref{102}, we obtain
\begin{equation} \label{103}\begin{array}{ll}
&\|\nabla
u\|_{L^2}^{2}+\displaystyle2\nu\displaystyle\int_{0}^{t}\|\Delta
u\|_{L^2}^{2}d\tau\\
\displaystyle &\hspace{0.3cm} \leq\|\nabla
u(0)\|_{L^2}^{2}+\left(\sup_{0\leq s\leq t}\|\nabla_{h}
u\|_{L^2}\right)\left(\displaystyle\int_{0}^{t}\|\nabla
u\|_{L^2}^{2}d\tau\right)^{\frac{1}{4}}\\&\ \ \ \ \
\times\left(\displaystyle\int_{0}^{t}\|\nabla_{h}\nabla
u\|_{L^2}^{2}d\tau\right)^{\frac{1}{2}}\left(\displaystyle\int_{0}^{t}\|\Delta
u\|_{L^2}^{2}d\tau\right)^{\frac{1}{4}}\\
\displaystyle &\hspace{0.3cm} \leq\|\nabla
u(0)\|_{L^2}^{2}+\displaystyle\left(\int_0^t\|u_3\|_{L^q}^{\frac{4p}{6-p}}\|\nabla
u\|_{{L^2}}^{2}d\tau\right)^{\frac{6-p}{4}}\displaystyle\left(\int_0^t\|\Delta u\|_{{L^2}}^{\frac{p-2}{2}}d\tau\right)^{\frac{p-1}{4}}\\
 &\hspace{0.5cm} + \displaystyle\|\nabla_{h}
u(0)\|_{L^2}^{2}\left(\displaystyle\int_{0}^{t}\|\Delta
u\|_{L^2}^{2}d\tau\right)^{\frac{1}{4}}\\
\displaystyle &\hspace{0.3cm} \leq\|\nabla
u(0)\|_{L^2}^{2}+\displaystyle\left(\int_0^t\|u_3\|_{L^q}^{\frac{4p}{6-p}}\|\nabla
u\|_{{L^2}}^{2}d\tau\right)^{\frac{6-p}{5-p}}\\
 &\hspace{0.5cm} + \displaystyle C\|\nabla_{h}
u(0)\|_{L^2}^{\frac{8}{3}}+ \nu\displaystyle\int_{0}^{t}\|\Delta
u\|_{L^2}^{2}d\tau.
\end{array}\end{equation}
We finally get
\begin{equation} \label{103}\begin{array}{ll}
&\|\nabla
u\|_{L^2}^{2}+\displaystyle\frac{3\nu}{2}\displaystyle\int_{0}^{t}\|\Delta
u\|_{L^2}^{2}d\tau\\ \displaystyle &\hspace{0.3cm} \leq\|\nabla
u(0)\|_{L^2}^{2}+\displaystyle C\|\nabla_{h}
u(0)\|_{L^2}^{\frac{8}{3}}+\displaystyle\left(\int_0^t\|u_3\|_{L^q}^{\frac{4p}{6-p}}\|\nabla
u\|_{{L^2}}^{2}d\tau\right)^{\frac{6-p}{5-p}}.
\end{array}\end{equation}
Now, combing \eqref{32} and \eqref{103}, again, by energy
inequality, and using H$\ddot{\mbox{o}}$lder's  and Young's
inequalities, we get
\begin{equation} \label{34}\begin{array}{ll}
\displaystyle&\|\nabla
u\|_{L^2}^{2}+\displaystyle\frac{3\nu}{2}\displaystyle\int_{0}^{t}\|\Delta
u\|_{L^2}^{2}d\tau\\ \displaystyle&\leq  C \left(\sup_{0\leq\tau\leq
t}\|u_3\|_{L^q}^2\right)^{\frac{2p}{5-p}}\displaystyle\left(\int_0^t\|\nabla
u\|_{{L^2}}^{2}d\tau\right)^{\frac{6-p}{5-p}}+\|\nabla
u(0)\|_{L^2}^{2}+\displaystyle\|\nabla_{h}
u(0)\|_{L^2}^{\frac{8}{3}}\vspace{1mm}\\
 \displaystyle&\leq  C \left(\sup_{0\leq\tau\leq
t}\|u_3\|_{L^q}^2\right)^{\frac{4q}{3q-10}}+\|\nabla
u(0)\|_{L^2}^{2}+\displaystyle\|\nabla_{h}
u(0)\|_{L^2}^{\frac{8}{3}}\vspace{1mm}\\
&\displaystyle\leq
C\left(\int_{0}^{t}\|\partial_{3}u_{3}\|_{\dot{B}_{\infty,\infty}^{\frac{2}{2-\beta}}}^{(1-\frac{2}{\beta})
 (\frac{\beta}{\beta-1})}\|\nabla
 u\|_{L^2}^{2}d\tau\right)^{\frac{4q(\beta-1)}{\beta(3q-10)}}
\left(\int_{0}^{t}\|\Delta u\|_{L^2}^{2}d\tau\right)^{\frac{4q}{\beta(3q-10)}}\vspace{1mm}\\
& \ \ \ \ \ \displaystyle +
C\|u_{3}(0)\|_{L^q}^{\frac{8q}{3q-10}}+\|\nabla
u(0)\|_{L^2}^{2}+\displaystyle\|\nabla_{h}
u(0)\|_{L^2}^{\frac{8}{3}}\vspace{1mm}\\
&\displaystyle\leq
C\left(\int_{0}^{t}\|\partial_{3}u_{3}\|_{\dot{B}_{\infty,\infty}^{\frac{2}{2-\beta}}}^{(1-\frac{2}{\beta})
 (\frac{\beta}{\beta-1})}\|\nabla
 u\|_{L^2}^{2}d\tau\right)^{\frac{4q(\beta-1)}{\beta(3q-10)-4q}}\displaystyle+C\|u_{3}(0)\|_{L^q}^{\frac{8q}{3q-10}}\vspace{1mm}\\
&\ \ \ \ \ +\|\nabla u(0)\|_{L^2}^{2}+\displaystyle\|\nabla_{h}
u(0)\|_{L^2}^{\frac{8}{3}}+\displaystyle\frac{\nu}{2}\displaystyle\int_0^t\|\Delta u\|_{{L^2}}^{2}d\tau.\\
\end{array}\end{equation}
In above inequality, we note that $\beta$ and $q$ satisfy the
additional condition
$$
\frac{4q}{\beta(3q-10)}<1.
$$
By H$\ddot{\mbox{o}}$lder's
and Young's inequalities, as well as  the energy inequality,  from \eqref{34} we have
\begin{equation} \label{35}\begin{array}{ll}
\displaystyle&\|\nabla
u\|_{L^2}^{2}+\displaystyle\nu\displaystyle\int_{0}^{t}\|\Delta
u\|_{L^2}^{2}d\tau\\&\displaystyle\leq
C\left(\int_{0}^{t}\|\partial_{3}u_{3}\|_{\dot{B}_{\infty,\infty}^{\frac{2}{2-\beta}}}^{\frac{4q(\beta-2)}{\beta(3q-10)-4q}}
\|\nabla u\|_{L^2}^{2}d\tau\right)\left(\int_{0}^{t}\|\nabla u\|_{L^2}^{2}d\tau\right)^{\frac{\beta(q+10)}{4q(\beta-1)}}\vspace{1mm}\\
& \ \ \ \ \ \displaystyle +
C\|u_{3}(0)\|_{L^q}^{\frac{8q}{3q-10}}+\|\nabla
u(0)\|_{L^2}^{2}+\displaystyle\|\nabla_{h}
u(0)\|_{L^2}^{\frac{8}{3}}\vspace{1mm}\\
&\displaystyle\leq
C\int_{0}^{t}\|\partial_{3}u_{3}\|_{\dot{B}_{\infty,\infty}^{\frac{2}{2-\beta}}}^{\frac{12(\beta-2)}{4\beta-37}}
\|\nabla u\|_{L^2}^{2}d\tau\vspace{1mm}\\
& \ \ \ \ \ \displaystyle +
C\|u_{3}(0)\|_{L^q}^{\frac{24\beta}{9\beta-25}}+\|\nabla
u(0)\|_{L^2}^{2}+\displaystyle\|\nabla_{h}
u(0)\|_{L^2}^{\frac{8}{3}}.\vspace{1mm}\\
\end{array}\end{equation}
Note that $\beta>\frac{37}{4}$, if we set $s=\frac{2}{\beta-2}$,
then we have
$$\frac{12(\beta-2)}{4\beta-37}=\frac{24}{8-29s}\ \mbox{with}\  0<
s<\frac{8}{29}.$$ By  using of  Gronwall's inequality, we obtain
\begin{equation}\label{36}\begin{array}{ll} \displaystyle &\displaystyle\|
\nabla u\|_{L^2}^{2}+\nu\int_{0}^{t}\|\Delta u\|_{L^2}^{2}d\tau\vspace{1mm}\\
&\displaystyle\leq \left(
C\|u_{3}(0)\|_{L^q}^{\frac{24\beta}{9\beta-25}}+\|\nabla
u(0)\|_{L^2}^{2}+\displaystyle\|\nabla_{h}
u(0)\|_{L^2}^{\frac{8}{3}}\right) \exp\left(
C\int_{0}^{t}\|\partial_{3}u_{3}\|_{\dot{B}_{\infty,\infty}^{-s}}^{\frac{24}{8-29s}}
d\tau\right),\\
\end{array}\end{equation}
by the condition   \eqref{qq}, the $ H^1$ norm of the strong
solution $u$ is bounded on the maximal interval of existence $(0,
T^{*})$. This completes the proof of Theorem \ref{t1.5}.

\section*{Acknowledgement}
The second author would like to thank Dr. Ting Zhang   for his
helpful suggestions.  This work is supported  by NSFC
11271322, 10931007, and Zhejiang NSF of China Z6100217.

}

\end{document}